\documentclass[12pt,leqno]{amsart}

\usepackage[utf8]{inputenc}   % LaTeX, comprends les accents !
\usepackage[T1]{fontenc}      % Police contenant les caract??res français
\usepackage{geometry}         % Définir les marges
\usepackage[francais]{babel}  % Placez ici une liste de langues, la
\usepackage{amssymb}
\usepackage{amsfonts}
\usepackage{amscd}
\usepackage[matrix]{xypic}
\usepackage{enumerate}
\voffset-1cm \markleft{\rm Poisson} 
\newtheorem{theorem}{Theorem}
\newtheorem{definition}[theorem]{Definition}
\newtheorem{corollary}[theorem]{Corollary}

\newtheorem{lemma}[theorem]{Lemma}

\newtheorem{proposition}[theorem]{Proposition}

\newtheorem{defi}{D\'{e}finition}

\newcommand{\K}{\mathbb K}

\newcommand{\N}{\mathbb N}
\newcommand{\R}{\mathbb R}

\newcommand{\KS}{{\mathbb K}[\Sigma_3]}
\newcommand{\ST}{\Sigma_3}

\newcommand{\ra}{\rightarrow}

\newcommand{\pf}{\noindent{\it Proof. }}
\newcommand{\ds}{\displaystyle}

\newcommand{\im}{{\rm Im }}

\newcommand\rk{\text{\rm rk}}

\newcommand{\Id}{\mathrm{id}}

\newcommand{\LE}{\mathcal{L}}
\newcommand{\vv}{a_1Id+a_2\tau_{12}+a_3\tau_{13}+a_4\tau_{23}+a_5c+a_6c^2}

\usepackage{graphicx}

\title{Deformation quantization of non associative algebras}
\author{ Elisabeth Remm}
\date{}
\address{}
\email{}

\begin{document}

\maketitle

\noindent{\bf Introduction}

In this work $\K$ is a field of characteristic $0$. By $\K$-algebra, we mean a $\K$-vector space $A$ with a bilinear map $\mu$ and we denote this algebra by $(A,\mu)$. We assume that $\mu$ satisfies a quadratic relation (but this can be extended to a $n$-ary relation) denoted $\mu \bullet \mu =0.$ For example, for the associative case, we have $\mu \bullet \mu = \mu \circ (\mu \otimes Id-Id \otimes \mu)$. The set of $n$-dimensional algebras satisfying a quadratic relation $\mu \bullet \mu =0$ is an algebraic variety $\mathcal{V}_n$ over $\K$ and the classical notion of formal deformation permits a description of  neighborhoods of any points of this variety (for a topology adapted to the structure of algebraic variety).  A naive definition of a formal deformation of a point $\mu \in \mathcal{V}_n$  is a formal series $\mu_t = \mu + \sum _{k \geq 1} t^k \varphi_k$, considered as a bilinear multiplication on the $\K[[t]]$-algebra $A[[t]]=A \otimes \K[[t]]$,  where the maps $\varphi_k$ are bilinear maps on $A$ which satisfy quadratic relations determinate by the formal identity
$$\mu_t \bullet \mu_t =0.$$
For example,  in degree $0$, $\mu \bullet \mu = 0$, in degree $1$, $\mu \bullet \varphi_1 + \varphi_1 \bullet \mu$ and so on.  Formal deformations are mainly used for the local study of $\mathcal{V}_n$. For example, 
a point of this variety such that all its deformations are isomorphic at this point is topologically rigid, that is, its orbit under the natural action of the linear group is open. But there are other applications of these deformations. For example, consider a formal deformation of a given point $\mu$. It determines new algebra laws that are related to the original law. In fact the linear term $\varphi_1$ of $\mu_t$ is also a multiplication on $A$ whose quadratic relation is a consequence of the relation of degree $1$, $\mu \bullet \varphi_1 + \varphi_1 \bullet \mu$. As a fundamental consequence, we can note the theory of deformation  quantization introduced in \cite{Li}. In a simplified way, if we consider a formal deformation  $\mu_t=\mu_0+t\varphi_1+\cdots$  of a commutative associative algebra $(A,\mu_0)$, the first term $\varphi_1$ is a cocycle for Hochschild cohomology associated with $(A,\mu_0)$. This first term is also an admissible Lie multiplication whose associated Lie bracket $\psi_1$ verifies the identity of Leibniz with the initial commutative associative law. Then this formal deformation determines naturally a Poisson algebra $(A,\psi_1,\mu_0)$ and $(A[[t]],\mu_t)$ is  the deformation quantization of the Poisson algebra  $(A,\mu_0,\psi_1)$. One of the goals of this work is to extend this construction for nonassociative algebras. As the world of nonassociative algebras is vast, we will focus on a class of nonassociative algebras whose quadratic definition relationship has symmetrical properties. These algebras have been previously studied in \cite{G.R.Nonass}.  

Also in order to see which algebras are associated with a given commutative algebra, we can use the process of polarization and depolarization introduced in \cite{RM} in the case of Poisson algebras. This process consists in looking, from a nonassociative multiplication, at the properties of symmetric and skew-symmetric bilinear applications that are attached to multiplication, which is a bilinear application, initially given.

\tableofcontents

\section{$v$-associative algebras}
Let  $\Sigma_3=\{Id,\tau_{12},\tau_{13},\tau_{23},c,c^2\}$ be the symmetric group of degree $3$ where $\tau_{ij}$ is the transposition between $i$ and $j$ and $c$ the cycle $(231)$. The multiplication $\sigma\sigma'$  corresponds to the composition $\sigma \circ \sigma'$. Let $\K[\Sigma_3]$ be the algebra group of $\Sigma_3$. It is provided with an associative algebra structure and with a $\Sigma_3$-module structure. The (left) action of $\ST$ on  $\KS$ is given by
$$(\sigma \in \ST, v=\sum a_i\sigma_i \in \KS) \ra \sum a_i \sigma  \sigma_i.$$ For any $v \in \KS$, $\mathcal{O}_l(v)=\{v,\tau_{12}v,\tau_{13}v,\tau_{23}v,cv,c^2v\}$ or more simply $\mathcal{O}(v)$ is the corresponding orbit and $F_v=\K[\mathcal{O}(v)]$ the $\K$-linear subspace of $\KS$ generated by $\mathcal{O}(v)$. It is also a $\ST$-submodule. 

\noindent Some notations:  We shall call the canonical basis of $\KS$ the family $\{Id,\tau_{12},\tau_{13},\tau_{23},c,c^2\}$ and $v=(a_1,a_2,a_3,a_4,a_5,a_6)$ the components of the vector $v$ in the canonical basis. We shall denote also by $M_v$ the matrix of the family $( v,\tau_{12}v,\tau_{13}v,\tau_{23}v,cv,c^2v)$ in the canonical basis:
$$
M_v=
\left(
\begin{array}{cccccc}
a_1 & a_2 & a_3 & a_4 & a_6 & a_5 \\
a_2	 & a_1 & a_5 & a_6 & a_4 & a_3 \\
a_3 & a_6 & a_1 & a_5 & a_2 & a_4 \\
a_4 & a_5 & a_6 & a_1 & a_3 & a_2 \\
a_5 & a_4 & a_2 & a_3 & a_1 & a_6 \\
a_6 & a_3 & a_4 & a_2 & a_5 & a_6 \\
\end{array}
\right)
$$

Let $(A,\mu)$ be a $\mathbb{K}$-algebra where $\mu$ is the multiplication  of $A.$
We note by $\mathcal{A}_{\mu}$ the associator of $\mu$ that is
$$\mathcal{A}_{\mu}= \mu \circ (Id \otimes \mu -\mu \otimes Id).
$$
Every $\sigma \in \Sigma_3$ defines a linear map 
$$\Phi_{\sigma} :  A^{\otimes ^3}  \rightarrow  A^{\otimes ^3}$$
given by
$$
\begin{array}{l}
\Phi_{\tau_{12}} (XYZ)=YXZ \\
\Phi_{\tau_{13}} (XYZ)=ZYX \\
\Phi_{\tau_{23}} (XYZ)=XZY \\
\Phi_{c} (XYZ)=YZX \\
\Phi_{c^2} (XYZ)=ZXY \\
\end{array}
$$
for any $X,Y,Z \in A$.
If $v=\sum_{\sigma \in \Sigma_3} a_{\sigma} \sigma \in \mathbb{K}\,[\Sigma_3]$, we shall put
$$\Phi_v=\sum_{\sigma \in \Sigma_3} a_{\sigma} \Phi_{\sigma}.$$
If $\sigma_1,\sigma_2 \in \Sigma_3$, then
$$\Phi_{\sigma_1} \circ \Phi_{\sigma _2}=\Phi_{\sigma_2\sigma_1}$$
and more generally, for any $v_1,v_2 \in \KS$:
$$\Phi_{v_1} \circ \Phi_{v _2}=\Phi_{v_2v_1}$$
\begin{definition}
An algebra $(A,\mu)$ is a $v$-algebra if there exists $v \in \mathbb{K}\, [\Sigma_3]$ such that
$$ \mathcal{A}_{\mu} \circ \Phi_v =0.
$$
\end{definition}
If $(A,\mu)$ is a $v$-algebra, then it is also a $v'$-algebra for any $v' \in \mathcal{O}(v)$. This implies that $(A,\mu)$ is a $v_1$-algebra for any $v_1 \in F_v$. But for any $v_1 \in F_v$ such that $\dim F_{v_1} < \dim F_v$, the $v_1$-associativity don't imply the $v$-associativity. For example, if $v=Id-\tau_{12} +c$, the vector $Id+\tau_{13} $ is in $F_v$. But $\dim F_{v_1}=3 < \dim F_v=4$ and the $(Id+\tau_{13})$-associativity doesn't imply the $(Id-\tau_{12} +v)$-associativity.

\medskip

There exists two particular vectors in $\KS$, denoted here $V_{Lad}$ and $V_{PA}$ (we shall see later the meaning of these notations), corresponding to the signum representation and the trivial representation:
$$V8{Lad}=\ds \sum_{\sigma \in \ST} \epsilon(\sigma)\sigma, \ \ V_{3Pa}=\ds \sum_{\sigma \in \ST} \sigma$$
where $\epsilon(\sigma)$ is the signup of the permutation $\sigma$. For each of these two vectors, the module $F_v$ is of dimension $1$ and they are the only vectors with this property.

\begin{proposition}
An algebra $(A,\mu)$ is Lie-admissible if and only if it is $V_{Lad}$-associative. An algebra $(A,\mu)$ is $3$-power-associative if and only if it is $V_{3Pa}$-associative.
\end{proposition}
The classes of Lie-admissible algebras and power-associative algebras have been introduced by Albert in \cite{Al}.  An algebra is called Lie-admissible if the skew bilinear map $\psi$ attached to $\mu$ is a Lie bracket. This is equivalent to write $\mathcal{A}_\mu \circ \Phi_{V_{Lad}}=0.$ An algebra  is said to be power-associative if the subalgebra generated by any element is associative. 
Over a field of characteristic $0$, an algebra is power-associative if and only if it satisfies $\mathcal{A}_{\mu}(x,x,x)=\mathcal{A}_{\mu}(x^2,x,x)=0$
for any $x \in A$. An algebra  is said to be $3$-power-associative  if and only if it satisfies $\mathcal{A}_{\mu}(x,x,x)=0$ for any $x \in A$. This last condition is equivalent to $\mathcal{A}_\mu \circ \Phi_{V_{3Pa}}=0.$ Let us note that the $3$-power-associativity implies
$\mathcal{A}_\mu(x^2,x,x)+\mathcal{A}_\mu(x,x^2,x)+\mathcal{A}_\mu(x,x,x^2)=0.$ But the condition $\mathcal{A}_\mu(x,x,x)=0$ implies $\mathcal{A}_\mu(x,x^2,x)=0$ and the  power-associativity is equivalent to $\mathcal{A}_\mu(x^2,x,x)+\mathcal{A}_\mu(x,x,x^2)=0.$ Then a sufficient condition for a $3$-power-associative to be power-associative is $\mathcal{A}_\mu \circ \Phi_{Id+\tau_{23}}= 0.$

\medskip

There is a third irreducible representation for the group $\ST$, the first two being associated with vectors $V_{Lad}$ and $V_{3Pa}$. This is the representation of degree $2$ 
 associated with the vector $v=Id+c+c^2$. Its orbit contains this vector and $\tau_{12}+\tau_{13}+\tau_{23}$. Then $\dim F_v=2.$ For example, $(A,\mu)$ with $\mu$ skew-symmetric is a $v$-algebra if and only if it is a Lie algebra. 

\section{Formal deformations of $v$-algebras}

Let $(A,\mu_0)$ be $v$-associative algebra, $v=\sum a_i\sigma_i \in \KS$.  Let $\varphi_1$ and $\varphi_2$ be two bilinear maps on $A$. We denote by 
$\varphi_1 \bullet_v \varphi_2$ the trilinear map on $A$:
$$\varphi_1 \bullet_v \varphi_2 (x_1,x_2,x_3)=\sum a_i(\varphi_1(x_{\sigma(1)},\varphi_2(x_{\sigma(2)},x_{\sigma(3)}))-\varphi_1(\varphi_2(x_{\sigma(1)},x_{\sigma(2)}),x_{\sigma(3)}))$$
for any $x_1,x_2,x_3 \in A$
that is
$$\varphi_1 \bullet_v \varphi_2=(\varphi_1 \circ (Id \otimes \varphi_2)-\varphi_1 \circ (\varphi_2 \otimes Id)) \circ \Phi_v.$$

Let $(A,\mu_0)$ be a $v$-algebra. A $v$-formal deformation of $(A,\mu_0)$ is given by a family of bilinear maps on $A$
$$\{\varphi_i : A \otimes A \rightarrow A , \  i \in \N\}$$
with $\varphi_0=\mu_0$ and
satisfying
\begin{equation}
\label{def}
\sum_{i+j=k, i,j \geq 0}\varphi_i \bullet_v \varphi_j=0, \ \ k\geq 0.
\end{equation}
If we denote by $\K[[t]]$ the algebra of formal series with one indeterminate $t$, this definition is equivalent  to consider on the space $A [[t]] $ of formal series with coefficients in $A$ a structure  of $\K[[t]]$-associative algebra such that the canonical map $A[[t]]/\K[[t]] \rightarrow A$ is an isomorphism of $v$-algebras. It is practical to write
$$\mu_t=\mu_0+t \varphi_1+t^2\varphi_2+ \cdots$$
Equation (\ref{def}) implies at the order $k=0$ that $\mu_0$ is $v$-associative. At the order $k=1$ we have
$$\varphi_0 \bullet_v \varphi_1 + \varphi_1 \bullet_v \varphi_0=\mu_0 \bullet_v \varphi_1 + \varphi_1 \bullet_v \mu_0=0.$$
To be consistent with conventional cohomological approaches to deformations, we will note by 
$\delta_{v,\mu_0}^2 \varphi$ the trilinear map
$$\delta_{v,\mu_0}^2 \varphi=\mu_0 \bullet_v \varphi_1 + \varphi_1 \bullet_v \mu_0.$$
In fact, we know that there is a cohomological complex which parametrizes formal deformations of quadratic algebras, and $\delta_{v,\mu_0}^2$ corresponds to the second coboundary operator. For example, if $v=Id$, then $(A,\mu_0)$ is associative and $\delta_{Id,\mu_0}\varphi$ is the coboundary operator associated with the Hochschild complex of $A$. It is also denoted in this case
$\delta_{H,\mu_0}$ and we have
$$\delta_{H,\mu_0}^2\varphi(x,y,z)=x\varphi(y,z)-\varphi(xy,z)+\varphi_1(x,yz)-\varphi(z,y)z$$
where, to simplify the notations, $xy$ means $\mu_0(x,y)$. Then
$$\delta_{v,\mu_0}^2 \varphi=\delta_{H,\mu_0}^2 \varphi \circ \Phi_v.$$Going back to Equation (\ref{def}), we have
\begin{equation}\label{eq}
\begin{array}{ll}
\medskip
 {\rm order} \  0 & \mathcal{A}(\mu_0) \circ \Phi_v =0,      \\
 \medskip
   {\rm order} \  1   &   \delta_{v,\mu_0 }^2\varphi_1 =0,\\
   \medskip
   {\rm order} \  2  & \varphi_1 \bullet_v \varphi_1 +  \delta_{v,\mu_0}^2 \varphi_2 =0.
\end{array}
\end{equation}
Therefore, there exists a vector $v_1 \in \KS$ such as $\varphi_1 \bullet_v \varphi_1 \circ \Phi_{v_1}$ if and only if 
$\delta_{v,\mu_0}^2 \varphi_2 \circ \Phi_{v_1}=0.$ But $\delta_{v,\mu_0}^2 \varphi_2 \circ \Phi_{v_1}=\delta_{H,\mu_0} \circ \Phi_{v_1v}=0.$ 
\begin{lemma}\label{1}
Let $(A,\mu_0)$ be a commutative algebra with $\mu_0 \neq 0$ and $\varphi$ a bilinear map on $A$.  If $\delta^2_{H,\mu_0}$ denote the Hochschild coboundary operator:
$$\delta^2_{H,\mu_0} \varphi (X,Y,Z)=X\varphi(Y,Z)-\varphi(XY,Z)+\varphi(X,YZ)-\varphi(X,Y)Z$$
where $XY$denotes the product $\mu_0(X,Y)$, then $\delta^2_{H,\mu_0}\varphi \circ \Phi_w= 0$ is equivalent to $w=V_{Lad}$.
\end{lemma}
\pf Let $w=a_1Id+a_2\tau_{12}+
+a_3\tau_{13}+a_4\tau_{23})+a_5c+a_6c^2$ be a vector of $\KS$. Since $\mu_0$ is commutative, then $\delta^2_{H,\mu_0}\varphi \circ \Phi_w$ corresponds to
$$
\begin{array}{l}
(a_1-a_5)X\varphi(Y,Z) +(a_4-a_5)X\varphi(Z,Y)+(a_2-a_4)Y\varphi(X,Z)+(a_5-a_6)Y\varphi(Z,X)\\
+(-a_1+a_6)Z\varphi(X,Y)+(-a_2+a_3)Z\varphi(Y,X)-(a_1+a_2)\varphi(XY,Z)-(a_3-a_5)\varphi(YZ,X)\\
-(a_4+a_6)\varphi(XZ,Y)+(a_1+a_4)\varphi(X,YZ)+(a_2+a_5)\varphi(Y,XZ)+(a_3+a_6)\varphi(Z,XY)
\end{array}$$
Then $\delta^2_{H,\mu_0}\varphi \circ \Phi_w=0$ if and only if $a_1=a_5=a_6=-a_2=-a_3=-a_4$ that is $w=a_1V_{Lad}$.

\medskip

Let us apply this lemme to study the equation
$$\varphi_1 \bullet_v \varphi_1 +  \delta_{v,\mu_0}^2 \varphi_2 =0.$$
For any $v_1 \in \KS$ we have
$$\varphi_1 \bullet_v \varphi_1\circ \Phi_{v_1} +  \delta_{v,\mu_0}^2 \varphi_2\Phi_{v_1}=\varphi_1 \bullet_{v_1v} \varphi_1 +  \delta_{v_1v,\mu_0}^2 \varphi_2 =0.$$
The vector $v$ given, there is an obvious solution to this equation corresponding to the case $v_1v=0.$  If $M_v$ is the matrix associated with $v$, then the equation $v_1v=0$ corresponds to the linear system $M_v V_1=0$ where $V_1$ is the column matrix of the vector $v_1$ and $v_1v=0$ if and only if $v_1 \in \ker  M_v$. 
But le rank of $M_v$ is the dimension of $\K[\mathcal{O}(v)]$. This rank is maximal if and only if $v \in \Sigma_3$. In this case $\mu_0$ is associative. In all the other cases, $\rk(M_v) < 6$ and $\dim \ker M_v \geq 1$.  For example, if $v=V_{Lad}$, then $\dim F_v=1$ and $v_1V_{Lad}=aV_{Lad}$ for any $v_1 \in \KS$. We deduce that, if $\mu_0$ is a commutative Lie admissible Lie algebra, then $\delta_{H,\mu_0}\varphi \circ \Phi_{v_1}=0$ for any $v_1$, implying $\delta_{H,\mu_0}\varphi =0$ for any bilinear map $\varphi$. We have similar result when $v=V_{3Pa}$, the second vector of rank $1$. In this case, $vV_{3Pa}=aV_{3Pa}$ implying that $\delta_{H,\mu_0}=0$ as soon as $\mu_0$ is a commutative $3$-power associative.

\noindent Let us assume now that $v_1 \notin \ker M_v$. From the previous lemma, $\delta^2_{H,\mu_0}\varphi \circ \Phi_{v_1v}=0$ if and only if $v_1v=kV_{Lad}$. The equation $v_1v=V_{Lad}$ is equivalent to the linear system $M_v V_1=V_{Lad}$ that is $V_{Lad} \in \im M_v$. 

\begin{proposition}\label{Mv}
Let $(A,\mu_0)$ be a $v$-commutative algebra and $\mu_t=\mu_0+t\varphi_1+t^2\varphi_2+ \cdots $ be a $v$-formal deformation of $\mu_0$. Then $\varphi_1$ is a Lie  admissible multiplication on $A$ if and only if $V_{Lad} \in \im M_v$.
\end{proposition}
\begin{corollary}\label{co}
Let $v=a_1Id+a_2\tau_{12}+a_3\tau_{13}+a_4\tau_{23}+a_5c+a_6c^2$ be a vector of $\KS$. If $\lambda=a_1-a_2-a_3-a_4+a_5+a_6 \neq 0$, then 
$V_{Lad} \in \im M_v$ and $\varphi_1$ is a Lie  admissible multiplication on $A$.
\end{corollary}
In fact,  if $\lambda=a_1-a_2-a_3-a_4+a_5+a_6 \neq 0$,  then $\lambda$ is a non null eigenvalue of $M_v$ and $V_{Lad}$ is an eigenvector attached to $\lambda$. In this case $V_{Lad} \in \im M_v$. 

Note that this condition on $\lambda$ is not necessary. An example is given by the Pre-Lie algebras, that is $v$-algebras with $v=Id-\tau_{23}$. For these algebras we have $\lambda=0$ and yet they are Lie-admissible. For the general case, in \cite{G.R.Nonass}, we determine vectors $v$ such that $V_{Lad} \in F_v$.

\section{Deformation quantization of the $v$-algebras}

Recall that  the rank of a vector $v \in \KS$ is the dimension of the vector space $F_v=\K[\mathcal{O}(v)]$. Then, if $v \neq 0$, then $1 \leq \rk(v) \leq 6.$ If $\rk(v)=6$ we have $F_v=\KS$ and $Id \in F_v$. In this case, any $v$-associative algebra is associative and we can assume that $v=Id$. Similarly, if $\rk(v)=1$, then $\dim F_v=1$ and it is  is a one dimensional invariant subspace of $\KS$. We have seen that, in this case, $v=V_{Lad}=Id-\tau_{12}-\tau_{13}-\tau_{23}+c+c^2$ or $v=V_{3Pa}=Id+\tau_{12}+\tau_{13}+\tau_{23}+c+c^2$. 

\subsection{Rank of $v$=6, $v=Id$: associative case}
The study of deformations of associative algebras was initiated by Gerstenhaber \cite{Ge}, and deformations quantizations by Bayen, Flato, Fronsdal, Lichnerowicz and Sternheimer in \cite{Li}. Here, in a first step,  we resume this study as part of the $v$-associative algebras

When $v=Id$, a $v$-algebra is an associative algebra. Let $\mu=\mu_0+t\varphi_1+t^2\varphi_2 + \cdots$ be an associative  formal deformation of a commutative associative multiplication $\mu_0$.  In this case $\delta_{H,\mu_0}^2=\delta_H^2$ and  Equations \ref{eq} writes
$$\begin{array}{ll}
 \medskip
   {\rm order} \  1   &   \delta_{H }^2\varphi_1 =0,\\
   \medskip
   {\rm order} \  2  & \varphi_1 \bullet \varphi_1 +  \delta_{H}^2 \varphi_2 =0.
\end{array}
$$
We have seen that for any bilinear map $\varphi$, any vector $v$ which satisfies $\delta_H^2 \varphi \circ \Phi_v=0$ is equal to $V_{Lad}$. 
So, for such a vector $v$ we have $\delta_H^2 \varphi_2 \circ \Phi_v=0$ what involves
$$\varphi_1 \bullet \varphi_1 \circ \Phi_v =0.$$
Since $v=V_{Lad}$, we deduce that $\varphi_1$ is Lie admissible. 

\begin{proposition} Let $\mu=\mu_0+t\varphi_1+t^2\varphi_2 + \cdots$ be an associative  formal deformation of a commutative associative multiplication $\mu_0$.Then $\varphi_1$ is a Lie admissible  multiplication. 
\end{proposition}
The bilinear map $\varphi_1$ satisfies also $\delta_{H }^2\varphi_1 =0$ and so $\delta_{H }^2\varphi_1\circ \Phi_v =0$ for any $v \in \KS$. Let’s determine this vector so that this relation involves a relation on the skew-bilinear map $\psi_1$ attached to $\varphi_1$, that is $\psi(x,y)=\varphi_1(x,y)-\varphi_1(y,x)$. If $v=a_1Id+a_2\tau_{12}+a_3\tau_{13}+a_4\tau_{23}+a_5c+a_6c^2$ 
with
$$a_5=a_1-a_3+a_4, \ \ a_6=a_1+a_2-a_3,$$
then, writing $xy$ for $\mu_0(x,y)$, we have  $xy=yx$ and
$$
\begin{array}{l}
\delta_{H }^2\varphi_1\circ \Phi_v =  a_1(\psi_1(xy,z)+\psi_1(xz,y)+\psi_1(zy,x))+a_2(y\psi_1(x,z)+z\psi_1(x,y)-\psi_1(xz,y) \\
 -\psi_1(xy,z))
 +a_3(\psi_1(xz,y)+x\psi_1(y,z)-z\psi_1(x,y))
 +a_4(-x\psi_1(y,z) -y\psi_1(x,z) 
 -\psi_1(xz,y)\\-\psi_1(zy,x))=0\\
\end{array}
$$
for any $a_1,a_2,a_3,a_4 \in \K$. This is equivalent to
$$
\left\{
\begin{array}{l}
\psi_1(xy,z)+\psi_1(xz,y)+\psi_1(zy,x)=0,\\
y\psi_1(x,z)+z\psi_1(x,y)-\psi_1(xz,y) -\psi_1(xy,z)=0,\\
\psi_1(xz,y)+x\psi_1(y,z)-z\psi_1(x,y)=0,\\
-x\psi_1(y,z) -y\psi_1(x,z) 
 -\psi_1(xz,y)-\psi_1(zy,x)=0.\\
 \end{array}
 \right.
 $$
The third identity is the Leibniz identity between the Lie bracket $\psi_1$ and the commutative associative multiplication $\mu_0$. Since the other identities are consequence of the Leibniz identity, we find again the classical result
\begin{proposition} If $\mu=\mu_0+t\varphi_1+t^2\varphi_2 +\cdots$ is an associative formal deformation of the commutative associative multiplication  $\mu_0$ on $A$, then $(A,\mu_0,\psi_1)$ is a Poisson algebra and the formal deformation $(A[[t]],\mu_t)$ is the deformation quantization of this Poisson algebra.
\end{proposition}

\medskip

In this proposition, we see that any associative deformation of $(A,\mu_0)$ gives a quantization. But, is there $v$-associative formal deformation of $\mu_0$, $v \in \KS$ but $v \notin \ST$ which defines a deformation quantization of a Poisson algebra $(A,\psi,\mu_0)$ for some Lie bracket $\psi$? Let $v$ be a vector of $\KS$. Since $\mu_0$ is associative, it is also $v$-associative. Then we can consider a $v$-deformation of $\mu_0$:
$$\mu=\mu_0+t\varphi_1 + t^2 \varphi_2 + \cdots$$
We deduce
\begin{equation}\label{eqv}
\begin{array}{ll}
\medskip
 {\rm order} \  0 & \mathcal{A}(\mu_0) \circ \Phi_v =0,      \\
 \medskip
   {\rm order} \  1   &   \delta_{v,\mu_0 }^2\varphi_1 =\delta_{H,\mu_0 }^2\varphi_1 \circ \Phi_v=0,\\
   \medskip
   {\rm order} \  2  & \varphi_1 \bullet_v \varphi_1 +  \delta_{H,\mu_0}^2 \varphi_2 \circ \Phi_v=0.
\end{array}
\end{equation}
We therefore want to determine all vectors $v\in \KS$ such that any $v$- deformation determines a Poisson algebra $(A,\mu_0,\psi)$. 
The skew bilinear map $\psi_1$ attached to $\varphi_1$ is  a Lie bracket, if and only if  $\varphi_1$ is Lie-admissible that is $ \varphi_1 \bullet_{V_{Lad}} \varphi_1=  \varphi_1 \bullet \varphi_1 \circ \Phi_{V_{Lad}}=0$.  But Equations (\ref{eqv}) imply that $\varphi_1 \bullet_v \varphi_1  =- \delta_{H,\mu_0}^2 \varphi_2 \circ \Phi_v$. From Proposition \ref{Mv}, there exists  $v_1 \in \KS$ be such that 
$$
\left\{
\begin{array}{l}
v_1v=v_{LAd},\\
 \delta_{H,\mu_0}^2 \varphi_2 \circ \Phi_{V_{Lad}}=\delta_{H,\mu_0}^2 \varphi_2 \circ \Phi_{v_1v} =0,\\
 \end{array}
 \right.
 $$
if and only if $V_{Lad} \in F_v$. 
Such vectors $v$ are described in \cite{G.R.Nonass}. Let us look now the consequences of the first equation $\delta_{H,\mu_0 }^2\varphi_1 \circ \Phi_v=0$. For this, let us consider as vector $v$, a vector of the following form $v=a_1Id+a_2\tau_{12}+a_3\tau_{13}+a_4\tau_{23}+(a_1-a_3+a_4)c+(a_1+a_2-a_3)c^2.$  We have
$$
\begin{array}{l}
\delta_{H }^2\varphi_1\circ \Phi_v =  a_1(\psi_1(xy,z)+\psi_1(xz,y)+\psi_1(zy,x))+a_2(y\psi_1(x,z)+z\psi_1(x,y)-\psi_1(xz,y) \\
 -\psi_1(xy,z))
 +a_3(\psi_1(xz,y)+x\psi_1(y,z)-z\psi_1(x,y))
 +a_4(-x\psi_1(y,z) -y\psi_1(x,z) 
 -\psi_1(xz,y)\\-\psi_1(zy,x)).\\
\end{array}
$$
To simplify the writing, we put
$$\LE(x,y,z)=\psi_1(xy,z)-x\psi_1(y,z)-y\psi_1(x,z)$$
that is the Leibniz identity between $\psi_1$ nd $\mu_0$. With these notations, the previous identity becomes
$$
\begin{array}{lll}
\delta_{H }^2\varphi_1\circ \Phi_v &= & -a_1(\LE(x,y,z)+\LE(y,z,x)+\LE(z,x,y))-a_2(\LE(x,z,y)+\LE(x,y,z))\\&& +a_3(\LE(x,z,y)) -a_4(\LE(x,z,y)+\LE(z,y,x))\\
&=& -(a_1+a_2)\LE(x,y,z)-(a_1+a_4)\LE(y,z,x)-(a_1+a_2-a_3+a_4)\LE(z,x,y).
\end{array}
$$
Then $\delta_{H }^2\varphi_1\circ \Phi_v=0$ implies $\LE=0$ if and only if
\begin{enumerate}
  \item $a_1+a_2 \neq 0$ and $a_1+a_4=a_2-a_3=0$,
  \item  or $a_1+a_4 \neq 0$ and $a_1+a_2=a_3-a_4=0$,
   \item or $a_3-a_4 \neq 0$ and $a_1+a_2=a_1+a_4=0.$
\end{enumerate}
For each case, the $v$ vector belongs to the orbit of $Id-\tau_{12}+c.$ Moreover $V_{Lad}$ belongs to $\K[\mathcal{O}(Id-\tau_{12}+c)]$.

\begin{proposition} Let $(A,\mu_0)$ be an associative commutative algebra and let be the vector $v=Id-\tau_{12}+c $. Then $(A,\mu_0)$ is $v$-associative and for any $v$-deformation $\mu_t=\mu_0+t\varphi_1+t^2\varphi_2 + \cdots$ of $\mu_0$, the algebra $(A,\mu_0,\psi_1)$ is a Poisson algebra where $\psi_1$ is the skew-symmetric application associated with $\varphi_1$.
\end{proposition}

\noindent{\bf Consequence.} Any deformation quantization constructed from a commutative associative algebra $(A,\mu_0)$ can be obtained by $v$-formal deformation, with $v=Id$ or $v=Id-\tau_{12}+c.$ In \cite{RMo}, the $v$-algebras with $v=Id-\tau_{12}+c$ are called weakly-associative algebras. An algebraic study of this class of algebras is presented in this paper.

\medskip

\subsection{Rank $1$: $v=V_{Lad}$ or $v=V_{3Pa}$}
Let $(A,\mu_0)$ be a commutative $V_{Lad}$-algebra and $\mu_t=\mu_0+ \sum t^i\varphi_i$ be a $V_{Lad}$-deformation of $\mu_0$. From Lemma \ref{1}, since $\mu_0$ is commutative, for any bilinear map $\varphi$ we have $\delta^2_{H,\mu_0}\varphi \circ \Phi_{V_{Lad}}=0$. This implies that , for any $i \geq 1$, $\delta^2_{V_{Lad},\mu_0} \varphi_i =0.$ In particular $\varphi_1$ is Lie admissible. 

\begin{proposition}
Let $(A,\mu_0)$ be a commutative Lie admissible algebra. Then for any any Lie-admissible formal deformation $\mu_t=\mu_0+t\varphi_1+t^2\varphi_2 +\cdots$, we have $\delta^2_{V_{Lad},\mu_0} \varphi_i =0$ and the algebra $(A,\varphi_1)$ is a Lie admissible.
\end{proposition}
Remark that, for any Lie admissible multiplication $\varphi$, then $\mu_t= \mu_0+t\varphi$ is a formal $V_{Lad}$-deformation, often called linear deformation.

\medskip

Let us consider now the $3$-power associative case. Let $(A,\mu_0)$ be a $3$-power associative commutative algebra.  For any bilinear application $\varphi$ on $A$ we have
$$
\begin{array}{l}
\delta_{V_{3PA }}^2\varphi (x,y,z)= -2(\psi(xy,z)+\psi(zx,y)+\psi(yz,x))
\end{array}
$$
where $\psi$ is the skew-bilinear map attached to $\varphi$. Let $\mu_t=\mu_0+\sum tî \varphi_i$ be a $V_{3Pa}$-formal deformation of $\mu_0$. In particular we have
$$\delta_{V_{3Pa }}^2\varphi_1=0,\ \ \ 
\varphi_1 \bullet_{V_{3Pa }} \varphi_1 +\delta^2_{V_{3Pa},\mu_0}\varphi_2=0.$$
But for any $v \in \KS$, we have $vV_{3Pa}=aV_{3Pa}$. We deduce that for any $v \in \KS$, 
$$(\varphi_1 \bullet_{V_{3Pa }} \varphi_1 +\delta^2_{V_{3Pa},\mu_0}\varphi_2)\circ \Phi_v=\varphi_1 \bullet \varphi_1 \circ \Phi_{vV_{3Pa}}+\delta^2_{H,\mu_0}\varphi_2\circ \Phi_{vV_{3Pa}}=a(\varphi_1 \bullet _{V_{3Pa}}\varphi_1 +\delta^2_{V_{3Pa},\mu_0}\varphi_2$$
and no identity on $\varphi_1 \bullet \varphi_1$ can be deduced.
\begin{proposition}
Let $(A,\mu_0)$ be a commutative (and then $3$ power associative) algebra. Then for any $3$-power associative formal deformation $\mu_t=\mu_0+t\varphi_1+t^2\varphi_2 +\cdots$ of $\mu_0$, $(A,\psi_1)$ is a skew-symmetric  algebra such that
$$\psi_1(xy,z)+\psi_1(yz,x)+\psi_1(zx,y)=0$$
for any $x,y,z \in A$ where $\psi_1$ is the skew-symmetric bilinear map attached to $\varphi_1$.
\end{proposition}
For example, this class of algebras contains the class of Poisson algebras. In fact if $(A,\mu_0,\psi)$ is a Poisson algebra, then
$$\LE(x,y,z)+\LE(y,z,x)+\LE(z,x,y)=\psi(xy,z)+\psi(yz,x)+\psi(zx,y)=0$$

\medskip

\subsection{Rank $2$ : $v=Id+c+c^2$, $v=Id-\tau_{12}+\tau_{13}-c^2$}
Since $F_v$ is a $\Sigma_3$-invariant vector space, it is a direct sum of irreducible vector spaces and the irreducible vectors spaces are of dimension $1$, that is $F_{V_{Lad}}$ and $F_{V_{3Pa}}$ or of dimension $2$, that is $F_{Id-\tau_{12}+\tau_{13}-c^2}.$

1. Assume that $F_v$ not irreducible. We can take $v=Id+c+c^2$
Any commutative multiplication $\mu_0$ satisfies $\mu_0 \bullet \mu_0 \circ \Phi_v
=0$ and $(A,\mu_0)$ is a $v$-algebra. We have, for any bilinear map $\varphi$ on $A$:
$$\delta^2_{v,\mu_0} \varphi (x,y,z)=\delta^2_{H,\mu_0} \varphi \circ \Phi_v(x,y,z) =\psi(x,yz)+\psi(y,zx)+\psi(z,xy)$$
for any $x,y,z \in A$, $\psi$ being the skew-symmetric map attached to $\varphi$. Let $\mu_t=\mu_0+\sum t^i\varphi_i$ be a $v$-formal deformation of $\mu_0$. Since we have
$$V_{Lad}\- v=3V_{Lad},$$
then $V_{Lad} \in F_v$ and $\varphi_1$ is a Lie admissible multiplication. Since $\delta^2_{v,\mu_0}\varphi_1=0$, the Lie bracket $\psi_1$ satisfies
$$\psi(x,yz)+\psi(y,zx)+\psi(z,xy)$$
for any $x,y,z \in A.$ 
\begin{proposition}
Let $v=Id+c+c^2$ be the vector of $\KS$ and  $(A,\mu_0)$  a commutative algebra. Then $(A,\mu_0$ is a $v$-algebra and for any $v$-formal deformation $\mu=\mu_0+t\varphi_1+t^2\varphi_2 +\cdots$ of $\mu_0$, $(A,\psi_1)$ is a Lie  algebra such that
\begin{equation}\label{r2}
\psi_1(xy,z)+\psi_1(yz,x)+\psi_1(zx,y)=0
\end{equation}
for any $x,y,z \in A$ where $\psi_1$ is the skew-symmetric bilinear map attached to $\varphi_1$.
\end{proposition}
We can call the triple $(A,\mu_0,\psi_1)$ a weakly Poisson algebra and $A[[t]],\mu_t)$ is  a deformation quantization of the weakly Poisson algebra $(A,\mu_0,\psi_1)$.  Any Poisson algebra is  weakly Poisson. An example of weak Poisson algebra but not Poisson is done considering the $2$-dimensional case: let $\{e_1,e_2\}$ be a basis of $A$ and
$$\psi_1(e_1,e_2)=e_2, \ \ \mu_0(e_1,e_1)=2\beta e_1, \ \mu_0(e_1,e_2)=\mu_0(e_2,e_1)=\alpha e_1+\beta e_2, \ \mu_0(e_2,e_2)=2\alpha e_2.$$
The algebra $(A,\psi_1,\mu_0$ is weakly Poisson and it is a Poisson algebra when $\alpha=0$ and $\beta=1.$ Recall that, in \cite{RM}, a Poisson algebra is represented by only one multiplication which satisfies a non associative identity, the two multiplications which appear in the definition of Poisson algebras being found again by the process of polarization-depolarization. This non associative multiplecation is called Poisson admissible.  If we apply this idea for weakly Poisson algebra, we find that the class of admissible weakly associative algebras corresponds to the $(Id+c+c^2)$-algebras.

\medskip

2. $F_v$ is irreducible. In this case none of the two vectors $V_{Lad}$ and $V_{3Pa}$ belongs to $F_v$. 

\subsection{Rank 3: $G_i$-algebras}
We will now focus mainly on $v$-algebras for which the linear term of a $v$-formal deformation is Lie admissible. We saw that this was equivalent to saying that $V_{Lad} \in F_v$. 
There is a particular class of $v$-algebras when $v$ is associated with a subgroup of $\ST$. The trivial subgroups, that is $G=\ST$ and $G=\{Id\}$ have been already studied. They concern the cases $\rk v= 1$ or $6$. When $G$ is the alternating subgroup, then $\rk v=2$ and this case was studied above. Then it remains the study concerning the subgroups 
$G_2=\{Id,\tau_{12}\},G_3=\{Id,\tau_{13}\}$ and $G_4=\{Id,\tau_{23}\}$.

\subsubsection{$G_2$-algebras: $v=Id-\tau_{12}$}. Let us begin by the case $v=Id-\tau_{12}$. A $v$-algebra is sometimes called a Pre-Lie algebra. Let $\mu_0$ be a commutative Pre-Lie algebra. For any bilinear map $\varphi$, we have
$$\delta^2_{v,\mu_0}(x,y,z)=x\varphi(y,z)-y\varphi(x,z)+\varphi(x,yz)-\varphi(y,xz)-z\psi(x,y)$$
for any $x,y,z \in A$. 
Let $\mu_t=\mu_0+\sum t^i\varphi_i$ be a $(Id-\tau_{12})$-formal deformation of $\mu_0$. Using the same notations as above, we have
$$
 \begin{array}{l}
   \delta^2_{v,\mu_0}\varphi_1=0,    \\
     \varphi_1 \bullet_v \varphi_1 +   \delta^2_{v,\mu_0}\varphi_1=0. 
\end{array}
$$
Since $V_{Lad} \in F_v$, the multiplication $\varphi_1$ is Lie admissible and, if $\psi_1$ is the skew-bilinear map attached to $\varphi_1$, $(A,\psi_1)$ is a Lie algebra.  Now, look for the identities concerning the Lie bracket $\psi_1$. Since 
$\delta^2_{v,\mu_0}\varphi_1=0$, we have $\delta^2_{v_1v,\mu_0}\varphi_1=0$ for any $v_1=\vv$. But 
$$v_1v=(a_1-a_2)(Id-\tau_{12})+(a_3-a_5)(\tau_{13}-c)+(a_4-a_6)(\tau_{23}-c^2)$$
implying that $v_1v \in F_v$ for any $v_1$.  Then the relation $\delta^2_{v_1v,\mu_0}\varphi_1=0$ is equivalent to $\delta^2_{v,\mu_0}\varphi_1=0$ and doen't give any new information.  Nevertheless, we can write again the expression $\delta^2_{v,\mu_0}\varphi_1=0$, using the Leibniz operator
$$\mathcal{L}(\mu_0,\varphi_1)(x,y,z)=\varphi_1(\mu_0(x,y),z)-\mu_0(x,\varphi_1(y,z))-\mu_0(\varphi_1(x,z),y)$$. We obtain
$$\delta^2_{v,\mu_0}\varphi_1=\mathcal{L}(\mu_0,\varphi_1) \circ \Phi_{Id-\tau_{12}}.$$
\begin{proposition}
Let $(A,\mu_0)$ be a commutative Pre-Lie algebra. Then any formal $(Id-\tau_{12})$-formal deformation $\mu=\mu_0+t\varphi_1 + \cdots$ determines an admissible Lie algebra $(A,\varphi_1)$ satisfying
$$\mathcal{L}(\varphi_1,\mu_0) \circ \Phi_{Id-\tau_{12}}=0.$$
\end{proposition}

\medskip

\subsubsection{$G_3$-algebra: $v=Id-\tau_{13}$}
A commutative $v$-algebra is also associative. Since $V_{Lad} \in F_v$, the linear term $\varphi_1$ of a $v$-formal deformation of $\mu_0$ is Lie admissible. The map $\varphi_1$ satisfies also
$$\delta_{v,H}\varphi_1(x,y,z)=x\rho_1(y,z)-z\rho_1(x,y)-\rho_1(xy,z)+\rho_1(x,yz)=0$$where $\rho_1$ is the symmetric map attached to $\varphi_1$.
 This is also written
$$\delta_{v,H}\varphi_1(x,y,z)=\rho_1(x,yz)-z\rho_1(x,y)-y\rho_1(x,z)-\rho_1(xy,z)+x\rho_1(y,z)+y\rho_1(x,z)=0$$
that is
$$\LE(\mu_0,\rho_1)(x,y,z)-\LE(\mu_0,\rho_1)(z,yx)=0.$$

\begin{proposition}
Let $(A,\mu_0)$ be a commutative  $(Id-\tau_{13})$-algebra. Then, if $\mu=\mu_0+t\varphi_1+t^2\varphi_2+\cdots$ is a $(Id-\tau_{13})$-formal deformation of $\mu_0$, then if $\psi_1$ and $\rho_1$ are respectively the skew-symmetric and symmetric bilinear maps attached to $\varphi_1$
\begin{enumerate}
\item $(A,\psi_1)$ is a Lie algebra,
\item The symmetric map $\rho_1$ satisfies
$$\LE(\mu_0,\rho_1)\circ \Phi_{Id-\tau_{13}}=0.$$
\end{enumerate}
\end{proposition}

\medskip

\noindent{Example.} Let $A=\mathcal{F}^1(\R,\R)$ be the commutative algebra of derivable real functions. Since it is associative, it is also $(Id-\tau_{13})$-associative for the natural product $\mu_0$. Let $\varphi_1$ the bilinear ap on $A$:
$$\varphi_1(f,g)=fg'.$$
Then $\psi_1(f,g)=fg'-f'g$ is a Lie braket and $\rho_1(f,g)=f'g+fg'=(fg)'$.  We have also
$$f\rho_1(g,h)-h\rho_1(f,g)-\rho_1(fg,h)+\rho_1(f,gh)=f(gh)'-h(fg)'-fgh'+f(gh)'=0.$$

\subsubsection{$G_4$-algebra: $v=Id-\tau_{23}$}
This case is similar to that of $G_2$-algebras.

\subsubsection{$v=Id+a\tau_{12}-\tau_{23}-ac^2$, \ $a \neq 1$}
From Proposition \ref{Mv}, $\varphi_1$ is Lie admissible. Moreover $\delta^2_{v,\mu_0}\varphi_1 = 0$ implies $\delta^2_{v_{Lad},\mu_0}\varphi_1 = 0$ 
because $v_{Lad}v=(2-2a)v_{Lad}$ and $\delta^2_{v_{Lad},\mu_0}\varphi_1 = 0$ implies 
$$x\psi_1(y,z)+y\psi_1(z,x)+z\psi(x,y)=0.$$
\begin{proposition}
The vector $v=Id+a\tau_{12}-\tau_{23}-ac^2$ is of rank $3$. If we assume $a \neq 1$, any $v$-formal deformation of a commutative $v$-algebra  is a deformation quantization of a weakly Poisson algebra $(A,\mu_0,\psi_1)$.
\end{proposition}

\subsubsection{$v=Id-\tau_{12}-2\tau_{23}+2c$}
This vector is also of rank $3$. It completes the classification given in \cite{G.R.Nonass}. For this vector we have a similar result to the previous case.

\subsection{Rank $4$}
The only really interesting case when $v$ is of rank $4$ is the one where $v=Id-\tau_{12}+c$  corresponding to algebras weakly associative algebra. This situation has been studied in later. Recall the result
\begin{proposition}
Any weakly associative deformation of a commutative algebra corresponds to a deformation quantization of a nonassociative Poisson algebra. Moreover,
any weakly associative deformation of a commutative associative algebra corresponds to a deformation quantization of a  Poisson algebra. 
\end{proposition}

\subsection{Rank $5$: $v=2Id-\tau_{12}-\tau_{13}-\tau_{23}+c$}
Let $\mu_t=\mu_O+\sum t^i\varpi_i$ be a $v$-formal deformation of the commutative $v$-algebra $(A,\mu_0)$. From Corollary \ref{co} $\varphi_1$ is Lie admissible. Since $\mu_0$ is commutative, we have
$$\begin{array}{l}
\delta^2_{v,\mu_0}\varphi(x,y,z)=x\varphi(y,z)+y\varphi(z,x)-2z\varphi(x,y)
-\varphi(xy,z)+\varphi(xz,y)+\varphi(x,yz)-\varphi(z,xy).
\end{array}
$$
If $v_1=(a_1,a_2,a_3,a_4,a_1+a_2-a_3,a_1+a_2-a_4)$, then $\delta^2_{v,\mu_0}\varphi\circ \Phi_{v_1}=0$ implies
$$
\begin{array}{l}
(a_2-2a_3+a_4)x\psi_1(z,y)+(-a_2-a_3+2a_4)y\psi_1(z,x)+(2a_2-a_3-a_4)z\psi_1(x,y) \\
+(a_4-a_3)\psi_1(xy,z)+(a_3-a_2)\psi_1(xz,y)+(a_2-a_4)\psi_1(yz,x)=0
\end{array}
$$
that is 
$$
\begin{array}{l}
a_2(-x\psi_1(y,z)-y\psi_1(z,x)+2z\psi_1(x,y) -\psi_1(xz,y)+\psi_1(yz,x))\\
+a_3(-2x\psi_1(y,z)-y\psi_1(z,x)-z\psi_1(x,y)-\psi_1(xy,z)+\psi_1(xz,y)\\
+a_4(x\psi_1(y,z)+2y\psi_1(z,x)-z\psi_1(x,y)+\psi_1(xy,z)-\psi_1(yz,x)=0.
\end{array}
$$
This is equivalent to the identity
$$-x\psi_1(y,z)-y\psi_1(z,x)+2z\psi_1(x,y) -\psi_1(xz,y)+\psi_1(yz,x)=0$$
that is 
$$\LE(\mu_0,\psi)(y,z,x)-\LE(\mu_0,\psi_1)(x,z,y)=0.$$
\begin{proposition}
Let $(A,\mu_0)$ be a $v$-algebra with $v=2Id-\tau_{12}-\tau_{13}-\tau_{23}+c.$ Any $v$-formal deformation $\mu_t=\mu_0+\sum t^i \varphi_i$ is a deformation quantization of a pseudo-Poisson algebra $(A,\mu_0,\psi_1)$ that is
\begin{enumerate}
\item $\mu_0$ is a commutative multiplication on $A$,
\item $\psi_1$ is a Lie bracket on $A$,
\item $\LE(\mu_0,\psi_1)(x,y,z)-\LE(\mu_0,\psi_1)(z,y,x)=0$.
\end{enumerate}
\end{proposition} 
%%%%%%%%%%%%%%%%%%%%%%%%%%%%%%%%%%%%%%%%%%%%%%%%%%%%%%%%%%%%%%%%%%%%%%%%
\section{A generalization : $\mathcal{V}$-$\mathcal{W}$-algebras}

This generalization has been introduced in \cite{G.R.Nonass}. If $(A,\mu)$ is a $\K$-algebra, we have noted by $\mathcal{A}(\mu)$ its associator. Let’s write this associator in the following form:
$$ \mathcal{A}_{\mu}=\mathcal{A}_{\mu}^R-\mathcal{A}_{\mu}^L
$$
where $\mathcal{A}_{\mu}^L=\mu \circ(Id \otimes \mu)$ and $\mathcal{A}_{\mu}^R=\mu\circ (\mu \otimes Id).$

Now, instead of considering action of $\Sigma_3$-permutation on the associator we can consider it independently on $\mathcal{A}_{\mu}^L$ and $\mathcal{A}_{\mu}^R$ which will induce different symmetries.

\subsection{Definition}

\begin{definition}
Let $v$ and $w$ two vectors of $\KS$. We say that the algebra $(A,\mu)$ is a $(v,w)$-algebra if we have
$$(\star)
\left\{
\begin{array}{l}
\mathcal{A}_{\mu}^L \circ \Phi_v=0 \\
\mathcal{A}_{\mu}^R \circ \Phi_w=0
\end{array}
\right.
\quad {\rm or} \quad (\star\star)
\left\{
\begin{array}{l}
\mathcal{A}_{\mu}^L \circ \Phi_v-\mathcal{A}_{\mu}^R \circ \Phi_w=0
\end{array}
\right.
$$
\end{definition}

We assume also that for every $v',w' \in \KS$
such that $v \in F_{v'}, \ w \in F_{w'}$ and $v' \notin F_v,\ w' \notin F_w$ we have $\mathcal{A}_{\mu}^L \circ \Phi_{v'} \neq 0$
or $\mathcal{A}_{\mu}^R \circ \Phi_{w'}\neq 0.$ In the study of $(v,w)$-algebras, the first study is to know if a $(v,w)$-algebra can be defined as a $v_1$-algebra. The simpler example is when $v=w$. In \cite{G.R.Nonass} we sudy the $(v,w)$-algebras which are Lie admissible algebras or Pre-Lie algebras. 

\medskip

\noindent{Example.} A Leibniz algebra satisfies the quadratic relation
$$
\mu(\mu(x,y),z)=\mu(x,\mu(y,z))+\mu(\mu(x,z),y)$$
 or
$$\mathcal{A}_{\mu}^L\circ \Phi_{Id-\tau_{23}}-\mathcal{A}_{\mu}^R \circ \Phi_{Id}= 0.$$
Then Leibniz algebras are $(Id-\tau_{23},Id)$-algebras. 
Symmetric Leibniz algebras are defined by a pair of quadratic relations. They correspond to
$$ \mathcal{A}_{\mu}^R \circ\Phi_{Id-\tau_{12}} -\mathcal{A}_{\mu}^L\circ \Phi_{Id}, \ \ \ \mathcal{A}_{\mu}^L \circ\Phi_{Id-\tau_{23}} -\mathcal{A}_{\mu}^R\circ \Phi_{Id},$$

\subsection{Formal deformations of $(v,w)$-algebras}
The notion of $(v,w)$-formal deformation of a $(v,w)$-algebra is similar to this notion for $v$-algebras.  Let $(A,\mu_0)$ be a $(v,w)$-algebra defined by a relation of type ($\star,\star$). Let $\mu_t= \mu_0 +\sum t^i\varphi_i$ be a formal deformation of $\mu_0$. We shall say that $\mu_t$ is a $(v,w)$-formal deformation of $\mu_0$ is $(A[[t]],\mu_t)$ is a $(v,w)$-algebra. To describe the relations between the $\varphi_i$, we need to introduce some notations:
$$\varphi_i \bullet_v^L \varphi _j = \varphi_i \circ (Id \otimes \varphi _j) \circ \Phi_v, \ \ \varphi_i \bullet_w^R \varphi _j = \varphi_i \circ (Id \otimes \varphi _j) \circ \Phi_v$$
$$\delta^{2,L}_{v,\mu_0}\varphi=\varphi \bullet^L_v \mu_0+ \mu_0 \bullet^L_{v} \varphi, \ \ \delta^{2,R}_{w,\mu_0}\varphi=\varphi \bullet^R_w \mu_0+ \mu_0 \bullet^R_{w} \varphi.$$
Thus, to say that $\mu_t$ is a $(v,w)$-formal deformation of $\mu_0$ implies in particular:
\begin{enumerate}
\item order $0$ : $\mu_0$ is a $(v,w)$-algebra,
\item order $1$ : $\delta^{2,L}_{v,\mu_0}\varphi_1+ \delta^{2,R}_{w,\mu_0}\varphi_1 =0$,
\item order $2$ : $\varphi_1 \bullet^L_v \varphi_1+\varphi_1 \bullet^R_w \varphi_1 + \delta^{2,L}_{v,\mu_0}\varphi_2
+ \delta^{2,R}_{w,\mu_0}\varphi_2=0$
\end{enumerate}

\subsection{Leibniz algebras}
Recall that $(A,\mu)$ is a Leibniz algebra if $\mu$ satisfies the quadratic relation
$$\mu(x,\mu(y, z)) = \mu(\mu(x, y ), z) + \mu(y, \mu(x, z))$$
for any $x,y,z \in A$ what is also written
$$\mathcal{A}(\mu)(x,y,z)= \mu(y, \mu(x, z)).$$

Let $(A,\mu_0)$ be a commutative Leibniz algebra. Writing $\mu_0(x,y)=xy$, we have
$$x(yz)-(xy)z-y(xz)=0.$$
For seach multiplication we have
$$\delta^2_{v,\mu_0}\varphi(x,y,z)=x\varphi(y,z)-y\varphi(x,z)-z\varphi(x,y)+\varphi(x,yz)-\varphi(y,xz)-\varphi(xy,z)$$
and 
$$\delta^2_{v,\mu_0}\varphi \circ \Phi_{v_1} =0 \Longrightarrow v_1=0.$$
So, if $\mu_t$ is a Leibniz-formal deformation (that is a $(Id-\tau_{12]},Id)$-formal deformation) of $\mu_0$, the relation
$$\varphi_1 \bullet^L_{Id-\tau_{12}} \varphi_1+\varphi_1 \bullet^R_{Id} \varphi_1 + \delta^{2,L}_{Id-\tau_{12},\mu_0}\varphi_2
+ \delta^{2,R}_{Id,\mu_0}\varphi_2=0$$
can not be reduced.
Let us consider now the relation 
$$\delta^{2,L}_{v,\mu_0}\varphi_1+ \delta^{2,R}_{w,\mu_0}\varphi_1 =0.$$
For any vector $v_1$ we have
$$(\delta^{2,L}_{v,\mu_0}\varphi_1+ \delta^{2,R}_{w,\mu_0}\varphi_1)\circ \Phi_{v_1} =0.$$
If $v_1=(a_1,a_2,a_3,-a_1,-a_2,-a3)$, this equation gives:
$$a_1x\psi_1(y,z)+a_2y\psi_1(x,z)-a_3z\psi_1(x,y)+(a_3-a_2)\psi_1(x,yz)-(a_1+a_3)\psi_1(y,xz) +(a_1+a_2)\psi_1(z,xy)$$
what is also written 
$$
\begin{array}{l}
a_1(x\psi_1(y,z)-\psi_1(y,xz)+\psi_1(z,xy))\\
+a_2(y\psi_1(x,z)-\psi_1(x,yz)+\psi_1(z,xy))\\
+a_3(-z\psi_1(x,y)+\psi_1(x,yz)-\psi_1(y,xz))= 0               à.
\end{array}
$$
and this is equivalent to the relation
$$x\psi_1(y,z)-\psi_1(y,xz)+\psi_1(z,xy)=0.$$
\begin{proposition}
Let $(A,\mu_0)$ be a commutative Leibniz algebra.  It is a $(Id- \tau_{12},Id)$-algebra and any $((Id- \tau_{12},Id)$-formal deformation of $\mu_0$ is a deformation quantization of a Pseudo-Poisson algebra $(A,\mu_0,\psi_1)$ that is
\begin{enumerate}
\item $\mu_0)$ be a commutative Leibniz multiplication,
\item $\psi_1$ is a skew-symmetric multiplication
\item we have the pseudo Leibniz relation
$$x\psi_1(y,z)-\psi_1(y,xz)-\psi_1(xy,z)=0$$
for any $x,y,z \in A$.
\end{enumerate}
\end{proposition}
\section{Polarization and depolarization}
Any multiplication $\mu : A 
 \otimes A \to A$ defined by a bilinear application can decomposed into the sum of a 
 commutative multiplication $\rho$ and an skew-symmetric one 
 $\psi$ via the {\em polarization \/} given by 
 \begin{equation}
\label{pol}
 \rho(x,y) =\frac{1}{2}(x y+y  x)\ 
 \mbox { and }\ 
 \psi(x,y) =\frac{1}{2}(x  y-y  x),\ 
 \mbox { for }\ x,y \in A. 
 \end{equation}
where $\mu(x,y)$ is denoted by $xy$.
 The inverse process of {\em depolarization\/} assembles a symmetric 
 multiplication $\rho$ with a skew-symmetric multiplication $\psi$ into 
 \begin{equation}
\label{depol}
 \mu(x,y) =\rho (x, y)+\psi(x,y),\ 
 \mbox { for }\ x,y \in V. 
 \end{equation}

 In the following section we give a well known example to illustrate the 
 (de)polarization trick. 
\subsection{Associative case} Assume that $(A,\mu)$ is an associative algebra.
 If we polarize the multiplication $\mu$
 we claim that decomposing $\mu(x,y)=\rho(x,y)+\psi(x,y)$, the associativity becomes
 equivalent to the following two axioms: 
 \begin{eqnarray} 
 \label{eq:3} 
 \psi(x,\rho(y ,z))&=& \rho(\psi(x,y), z) +\rho(y,\psi(x,z)), 
 \\ 
 \label{eq:1} 
 \psi(y,\psi(x,z)) &=&-\mathcal{A}(\rho) (x,y,z). 
 \end{eqnarray}
 To verify this, observe that associativity is equivalent to
$$R(\psi,\rho)=(\psi+\rho)\circ(Id\otimes \psi-\psi\otimes Id)+(\psi+\rho)(Id \otimes  \rho- \rho \otimes Id)=0.$$ 
Then $R(\psi,\rho) \circ \Phi_v=0$ for any $v \in \KS$. In particular
$$R(\psi,\rho) \circ \Phi_{V_{Lad}})=\psi \circ (\psi \otimes Id) \circ \Phi_{Id+c+c^2}= 0$$
and $\psi$ is a Lie bracket. The relation (\ref{eq:3}) shows that $(A,\rho,\psi)$ is a non associative Poisson algebra. It is a Poisson algebra if and only if the Lie bracket $\psi$ is $2$-step nilpotent.

Although the associative case is well known, we will try to put in place a systematic method to solve this case in order to extrapolate it to the other identities that interest us. Let be $\vv \in \KS$ and let us consider the identity $R(\psi,\rho) \circ \Phi_v=0.$ By grouping the terms $\psi(\Id\otimes \psi)$, $\rho(\Id\otimes \rho)$, $\rho(\Id\otimes \psi)$ and finally $\psi(\Id\otimes \rho)$, the coefficients of each of these terms are given by the matricial equation
$$N_v \left(
\begin{array}{c}
a_1\\
a_2\\
a_3\\
a_4\\
a_5\\
a_6
\end{array}
\right)$$
where $N_v$ is the transpose of the matrix
$$ \left(
\begin{array}{cccccccccccc}
1&0&1&1&0&-1&1&0&-1&1&0&1\\
0&1&-1&0&1&-1&0&1&1&0&1&1\\
-1&0&-1&-1&0&1&1&0&-1&1&0&1\\
-1&1&0&1&-1&0&-1&-1&0&1&1&0\\
1&-1&0&-1&1&0&-1&-1&0&1&1&0\\
0&-1&1&0&-1&1&0&1&1&0&1&1
\end{array}
\right)
$$
The rank of $M_v$ is $6$. Let us search the vectors of this space associated with minimal relations, that is to say with a maximum of 0 among these components. We obtain the independent vectors
$$\left\{
\begin{array}{l}
(0,0,0,0,0,0,0,-1,-1,1,0,0)\\
(0,0,0,0,0,0,-1,0,1,0,1,0)\\
(0,0,0,0,0,0,1,1,0,0,0,1)\\
\end{array}
\right.
$$
correspond to the vector $v=(1,-1,1,1,1,-1)$ and $\tau_{12}v,\tau_{13}v$ and the relation
$$\psi(x_1,\rho(x_2,x_3))-\rho(x_2,\psi(x_1,x_3))-\rho(x_3,\psi(x_1,x_2))=0.$$
This relation can be written
$\mathcal{L}(\psi,\rho)=0$ where $\psi$ is a Lie bracket and $\rho$ a commutative (nonassociative) multiplication. 
Similarly, we have the three  independent vectors
$$\left\{
\begin{array}{l}
(0,1,0,1,0,-1,0,0,0,0,0,0)\\
(1,0,0,0,1,-1,0,0,0,0,0,0)\\
(-1,1,1,0,0,0,0,0,0,0,0,0)\\
\end{array}
\right.
$$
which correspond to $v=(1,1,-1,1,-1,-1)$ and $\tau_{12}v,\tau_{13}v$ and to the relation
$$\psi(x_1,\psi(x_2,x_3))-\rho(x_3,\rho(x_1,x_2))+\rho(\rho(x_3,x_1),x_2).$$
\begin{proposition} Any associative algebra is associated with the principle of polarization 
depolarization to a bialgebra $(A,\psi,\rho)$ where $\psi$ is a Lie bracket, $\rho$ a commutative multiplication satisfying
\begin{enumerate}
\item $\mathcal{L}(\psi,\rho)=0$
\item $\mathcal{A}(\psi)=\mathcal{A}(\rho).$
\end{enumerate}
\end{proposition}
In particular, if $\rho$ is associative, then $(A,\psi,\rho)$ is a Poisson algebra whose the Poisson bracket is $2$-step nilpotent.

\subsection{Lie-admissible case}Recall that a nonassociative algebra is Lie admissible if its attached skew-symmetric bilinear map 
is a Lie bracket. In this case the principle of polarization says nothing else.

\subsection{Pre-Lie algebras}
A nonassociative algebra $(A,\mu)$ is a Pre-Lie algebra if its associator satisfies
$$\mathcal{A}_\mu(x,y,z)-\mathcal{A}_\mu(y,x,z)=0.$$
Since $V8{Lad} \in \K[\mathcal{O}(Id-\tau_{12})]$, such algebras are Lie admissible. If we polarize the multiplication $\mu$, we obtain
$$
\begin{array}{l}
\psi(x_3,\psi(x_1,x_2)) +\psi(x_1,\rho(x_2,x_3))-\psi(x_,\rho(x_1,x_3))+\rho(x_1,\psi(x_2,x_3))+\rho(x_2,\psi(x_3,x_1))-
\\
-2\rho(x_3,\psi(x_1,x_2))+\rho(x_1,\rho(x_2,x_3))-\rho(x_2,\rho(x_1,x_3))=0.
\end{array}
$$
Following the same procedure as in the previous case, the rank of the matrix $M_v$  is $3$ and  vectors of a basis of  the  image of $M_v$ are in the orbit. This implies that there is no relation other than that given in polarization.

\subsection{$\tau_{13}$-algebras}
These are the algebras $(A,\mu)$ defined by the quadratic relation
$$\mathcal{A}_{mu}(x_1,x_2,x_3)-\mathcal{A}_{mu}(x_3,x_2,x_1)=0.$$This relation is equivalent to
$$\psi(x_1,\psi(x_2,x_3))-\mathcal{A}_{\rho}(x_3,x_1,x_2)=0$$
and $\psi$ is a Lie bracket. This relation is minimal. 

\subsection{$(Id+c+c^2)$-algebras} 
In this case we have
$$\mathcal{A}_{\mu}(x_1,x_2,x_3)+\mathcal{A}_{\mu}(x_2,x_3,x_1)+\mathcal{A}_{\mu}(x_3,x_1,x_2)=0.$$
In this case $\psi$ is a Lie bracket and the polarization principe gives
$$\psi(x_1,\rho(x_2,x_3))+\psi(x_2,\rho(x_3,x_1))+\psi(x_3,\rho(x_1,x_2))=0.$$
It is a similar identity as that obtained by deformation. Let us note also that a $(Id+c+c^2)$-algebra is Lie admissible and $3$-power associative. 

\subsection{Weakly associative algebras} 
This class of nonassociative algebras has been studied in \cite{MoR} in view to extend the notion of quantification by deformation of associative commutative algebras. Recall that a nonassociative algebra $(A,\mu)$ is weakly associative if we have
$$\mathcal{A}_{\mu}(x_1,x_2,x_3)+\mathcal{A}_{\mu}(x_2,x_3,x_1)-\mathcal{A}_{\mu}(x_2,x_1,x_3)=0.$$
With reference to \cite{MoR}, it can be said that this identity is equivalent to:
\begin{enumerate}
\item $\psi$ is a Lie bracket,
\item $\rho$ is a commutative multiplication satisfying 
$$\psi(x_1,\rho(x_2,x_3)-\rho(x_2,\psi(x_1,x_3))-\rho(x_3,\psi(x_1,x_2)).$$
\end{enumerate}
In other words, $(\psi,\rho)$ satisfy the Leibniz identity : $\mathcal{L}(\psi,\rho)=0.$
\begin{proposition}
Let $(A,\mu)$ a weakly associative algebra and $(A,\rho,\psi)$ its depolarized. Then 
\begin{enumerate}
\item $(A,\psi)$ is a Lie algebra,
\item $\rho$ is a commutative multiplication satisfying 
$$\psi(x_1,\rho(x_2,x_3)-\rho(x_2,\psi(x_1,x_3))-\rho(x_3,\psi(x_1,x_2))$$
that is $(A,\rho,\psi)$ is a nonassociative Poisson algebra.
\end{enumerate}
\end{proposition}
\noindent

\noindent{Remark.} If we refer to \cite{G.R.Nonass}, weakly associativity corresponds to a point of the family of nonassociative algebras corresponding to the identity
$$\mathcal{C}_\mu(\alpha)(x,y,z)=2\mathcal{A}_\mu (x,y,z)+(1+\alpha)\mathcal{A}_\mu (y,x,z)+\mathcal{A}_\mu (z,y,x)+\mathcal{A}_\mu (y,z,x)+(1-\alpha)\mathcal{A}_\mu (z,x,y) =0$$
with $\alpha =-1/2.$

\medskip

If we consider the vector $v_1=\frac{1}{3}Id-\tau_{12}+\frac{7}{12}\tau_{13}+\frac{1}{4}c^2$, then the polarization of $\mathcal{C}_{\mu}(\alpha) \circ \Phi_{v_1}$ gives the relation
$$\mathcal{L}(\psi,\rho)(x_1,x_2,x_3)-\gamma(\psi(x_1,\psi(x_2,x_3))-2\psi(x_3,\psi(x_1,x_2))$$
with $\gamma=\frac{2}{3}(2\alpha-1)$. Since $v_1$ is  inversible in the algebra $\KS$, this relation is equivalent to $\mathcal{C}_\mu(\alpha)=0.$ 
\begin{proposition}
Let $v=2Id+(1+\alpha)\tau_{12}+\tau_{13}+c+(1-\alpha)c^2$ with $\alpha \neq 1$. Then any $v$-algebra is Lie admissible and $3$-power associative. The relation $\mathcal{A}_\mu \circ \Phi_v =0$ is equivalent to
$$\mathcal{L}(\psi,\rho)(x_1,x_2,x_3)-\gamma(\psi(x_1,\psi(x_2,x_3))-2\psi(x_3,\psi(x_1,x_2))$$
with $\gamma=\frac{2}{3}(2\alpha-1)$. In particular, if $\alpha=\frac{1}{2}$, then  $(A,\mu)$ is weakly associative and  we have in this case
$$\mathcal{L}(\psi,\rho)=0$$
that is $(A,\psi,\rho)$ is a nonassociative Poisson algebra.
\end{proposition}

\subsection{Leibniz algebras}
Recall that a Leibniz algebra is a quadratic algebra whose multiplication $\mu_ {x,y}=xy$ satisfies the identity
$$x(yz)-(xy)z-y(xz)=0.$$
Let $(\rho,\psi)$ be the pair of bilinear map given by the polarization of $\mu$.
They satisfy the identity:
$$\begin{array}{l}
R(\rho,\psi)(x,y,z)=\rho(x,\rho(y,z))-\rho(z,\rho(x,y))-\rho(y,\rho(x,z))+\rho(x,\psi(y,z))-\rho(z,\psi(x,y))\\-\rho(y,\psi(x,z))
+\psi(x,\rho(y,z))+\psi(z,\rho(x,y))+\psi(y,\rho(x,z))+\psi(x,\psi(y,z))\\+\psi(z,\psi(x,y))+\psi(y,\psi(x,z))
=0.
\end{array}
$$
Let $v=(a_1,a_2,a_3,a_4,a_5,a_6)$ be in $\KS$ where the $a_i$ are the components in the canonical basis. The matrix of the linear system $R(\rho,\psi) \circ \Phi_v$ is
$$N_v=
\left(
\begin{array}{rrrrrr}
1&-1&-1&1&-1&-1\\
-1&1&-1&-1&1&-1\\
-1&-1&1&-1&-1&1\\
1&-1&1&-1&-1&1\\
-1&1&1&-1&-1&1\\
-1&1&-1&-1&1&1\\
1&-1&1&1&1&-1\\
-1&1&-1&1&1&1\\
1&1&1&-1&-1&1\\
1&-1&-1&-1&1&1\\
1&-1&-1&-1&1&1\\
1&-1&-1&-1&1&1\\
\end{array}
\right)
$$
Its rank is equal to $6$. Let us determine a "`good"' basis of $\im N_v$.
The vectors column $C_3+C_5,C_4+C_6,C_1+C_2$ give the vectors
$$\left\{
\begin{array}{l}
u_1=\-^t(-1,0,0,0,0,0,1,0,0,0,0,0)\\
u_2=\-^t(0,-1,0,0,0,0,0,1,0,0,0,0)\\
u_3=\-^t(0,0,-1,0,0,0,0,0,1,0,0,0)
\end{array}
\right.
$$
The corresponding identities are
$$R(\rho,\psi)\circ \Phi_{w_i}=0$$
with
$$w_1=Id+\tau_{12}, \ w_2=\tau_{13}+c, \ w_3=\tau_{23}+c^2.$$
Sinc $w_2,w_3 \in \mathcal{O}(w_1)$, these identities are equivalent to $R(\rho,\psi)\circ \Phi_{w_1}=0$, that is $$\rho(x,\rho(y,z))=\psi(x,\rho(y,z))$$
for any $x,y,z \in A$.
For the other three vectors of the $\im N_v$, let’s take
$u_4=C_1-C_2+C_3-C_4+C_5-C_6$ which corresponds to the vector of $\KS$:
$$w_4=Id-\tau_{12}+\tau_{13}-\tau_{23}+c-c^2$$
We consider also the other vectors of $\mathcal{O}(w_4)$, that is $w_5=\tau_{13}w_4,w_6=\tau_{23}w_4.$ So we obtain a basis of $\im N_v$ and a second relation
$R(\rho,\psi)\circ \Phi_{w_4}=0$.
\begin{proposition}
Let $(A,\mu)$ a Leibniz algebra. It is associated, from the polarization - depolarization principle with a bialgebra $(A,\rho,\psi)$ where $\rho$ is a commutative multiplication, $\psi$ a skew-symmetric multiplication satisfying
\begin{enumerate} \item $\rho(x,\rho(y,z))=\psi(x,\rho (y,z))$,
\item $\begin{array}{l}
\rho(x,\psi(y,z))-\rho(y,\psi(x,z))-\rho(z,\psi(x,y))+2\psi(x,\rho(y,z))
-2\psi(y,\rho(x,z)\\+J(\psi)(x,y,z)=0
\end{array}$
\end{enumerate}
where $J(\psi)$ is the Jacobiator ( $\psi \circ  (Id \otimes \psi ) \circ \Phi_{Id+c+c^2})$, for any $x,y,z \in A$.
\end{proposition}

\medskip

\noindent{\bf Remark: Case of symmetric Leibniz algebras}.  Recall that such algebras correspond to the two identities:
$$
\left\{
\begin{array}{l}
x(yz)-(xy)z-y(xz)=0,\\
(xy)z-x(yz)-(xz)y=0
\end{array}
\right.
$$
This pair of relations is equivalent to
$$\mathcal{A}(\mu)(x,y,z)-y(xz)=0, \ \mathcal{A}(\mu)(y,z,x)+(yx)z=0$$
that implies
$$\mathcal{A}(\mu)(x,y,z)+ \mathcal{A}(\mu)(y,z,x)- \mathcal{A}(\mu)(y,x,z=0).$$
\begin{proposition}
Any symmetric Leibniz algebra is weakly associative. In particular $(A,\rho,\psi)$ is a nonassociative Poisson algebra. 
\end{proposition}

\end{document}